# On the machine swing dynamics: a perspective

Prashant G. Medewar† and Shambhu N. Sharma‡

*Department of Electrical Engineering,*

*Sardar Vallabhbhai National Institute of Technology, Surat, Gujarat, India,*

E-mails: †medewar.prashant@gmail.com and ‡snsvolterra@gmail.com

**Abstract:** A formal approach to rephrase nonlinear filtering of stochastic differential equations is the Kushner setting in applied mathematics and dynamical systems. Thanks to the ability of the Carleman linearization, the 'nonlinear' stochastic differential equation can be equivalently expressed as a finite system of 'bilinear' stochastic differential equations with the augmented state under the finite closure. Interestingly, the novelty of this paper is to embed the Carleman linearization into a stochastic evolution of the Markov process. To illustrate the Carleman linearization of the Markov process, this paper embeds the Carleman linearization into a nonlinear swing stochastic differential equation. Furthermore, we achieve the nonlinear swing equation filtering in the Carleman setting. Filtering in the Carleman setting has simplified algorithmic procedure. The concerning augmented state accounts for the nonlinearity as well as stochasticity. We show that filtering of the nonlinear stochastic swing equation in the Carleman framework is more refined as well as sharper in contrast to benchmark nonlinear EKF. This paper suggests the usefulness of the Carleman embedding into the stochastic differential equation to filter the concerning nonlinear stochastic differential system. This paper will be of interest to nonlinear stochastic dynamists exploring and unfolding linearization embedding techniques to their research.

# 1. Introduction

The problem of achieving the estimation and control of nonlinear dynamical systems poses mathematical intractability attributed to the infinite dimensionality of system dynamics filtering (Sharma *et al.* [1]). To circumvent the curse of dimensionality, linearization has proven useful. The Carleman linearization has formal and systematic constructions with its ability to preserve the nonlinearity as well as offer the simplified dynamical structure. Here, we explain partially the beauty, power and universality of the Carleman linearization in the dynamics and control. The Riccati equation has found its connections in dynamical and controls revealing different matrix interpretations. Interestingly, the Carleman linearization has connections with the Riccati equation. The Carleman embedding into the Riccati equation, a nonlinear ordinary differential equation, results into a finite system of linear equations under the finite closure (Bellman [2]). Despite the long-standing history of the Carleman linearization, greater conceptual depth as well as its universality for the nonlinear dynamical system, an extension of the Carleman embedding to the nonlinear stochastic differential system, is relatively very less researched in theory as well as its applications to practical problems. The Carleman linearization has found its application in the nonlinear stochastic networks (Bhatt and Sharma [3]).

The power system dynamics have nonlinearity coupled with stochasticity as well. The nonlinear swing equation is the cornerstone in the power system for the stability of the equilibrium point under the perturbations. That is well studied in the literature (Lu *et al.* [4]). It plays dominant roles in power systems, e.g. rotor angle stability in the Single Machine Infinite Bus (SMIB) system as well as coupled multiple machine systems, low-frequency oscillation, and sub synchronous oscillation. This equation is non-trivially nonlinear. Small local disturbances in power systems can lead to consequences influencing the whole system. A survey on power system blackouts and causes is given in Alhelou *et al.* ([5]). To gain insights into swing recording equipment's data as well as extracting states and parameters, algorithmic mathematics of data needs refinements. That can be achieved using two ideas. The first is to embed the stochasticity into the system dynamics. Since the dynamics is nonlinear and stochastic, linearization of the nonlinear dynamical system is the cornerstone to make the concerning computations tractable. Then the estimation matrix algebra and nonlinear filtering stochastic differential equations are the major ingredients to sketch the algorithmic procedures for extracting states and parameters from the swing recording equipment's data.

The objectives of the paper are two-fold. The first is to achieve 'the Carleman embedding into the nonlinear stochastic evolution of the Markov process'. The Markov process has independent increments as well as practical stochastic systems have generally the Markovian character. The second is to filter nonlinear machine swing dynamics in the Carleman setting. The first is achieved via a unification of the generating function, Carleman linearization, stochastic differential rules and the finite closure. The machine swing stochastic dynamics is the Markovian. Concerning the second objective, in this paper, we recast the finite-dimensional nonlinear swing equation into a finite system of bilinear stochastic differential equations via the Carleman embedding. After using the notion of matrix partitioning as well as stochastic differential equations associated with the Carleman linearized system dynamics, we sketch the swing equation filtering in the Carleman framework. To contrast the utility of the swing equation filtering in the Carleman setting, we exploit the celebrated extended Kalman filtering. Numerical simulations for the two sets of system parameters as well as initial conditions unfold the superiority of the 'swing equation filtering' in the Carleman setting in contrast to the swing equation filtering in the EKF setting. Perhaps, this is the first paper in two senses. (i) extending the idea of the Carleman linearization to nonlinear stochastic differential systems without the loss of generality in the system nonlinearity and nonlinear process noise diffusion coefficient. (ii) accomplishing its application via achieving machine equation filtering in the Carleman setting.

The rest of the paper is organized as follows: Section 2 develops a new theoretical result that transforms the nonlinear SDE into a bilinear SDE using the idea of the Carleman linearization. Section 3 provides a discussion embedding the Carleman linearization into the nonlinear stochastic SMIB system, i.e. a nonlinear noisy machine swing dynamics. In addition, the noisy machine swing dynamics filtering in Carleman setting as well as in EKF setting. Section 4 presents quite extensive simulation results with two sets of parameters for the noisy machine swing dynamics and it shows the superiority of the proposed Carleman setting of the algorithm. Concluding remarks are presents in Section 5.

# 2. Carleman linearization of a Markov process

Here, we explain succinctly the Carleman linearization of the stochastic evolution of a Markov process. For the simplified analysis as well as gaining insights, consider a nonlinear scalar Itô stochastic differential equation and then embed the Carleman linearization. The Carleman embedding results in an infinite system of bilinear stochastic differential equations. Finite closure circumvents the curse of dimensionality. As a consequence of this, we arrive at a finite system of bilinear stochastic differential equations. Notably, the finite closure preserves the Markovity into the Carleman linearized stochastic differential equation as well. The Carleman linearization brings bilinearity as well as preserves the nonlinearity via associating the 'nonlinearity stochastic evolution' with the 'stochastic state evolution'. Secondly, for the scalar stochastic process, the dimension of the Carleman linearized state vector becomes the order of the Carleman linearization order $N.$ For the vector case, the dimension becomes $\sum_{1\le r\le N}\binom{n+r-1}{r}$, where $n$ is the dimension of the state vector and $r$ is the Kronecker power associated with the state that runs over $1$ to $N.$

The cubic nonlinearity, which has received attention in the literature, e.g. Jing and Lang [6], is relatively very general in contrast to the square. Thus, the Carleman framework of the nonlinear stochastic differential equation preserving the cubic nonlinearity is sketched in the following Theorem-proof format.

**Theorem**: Consider a Markovian stochastic evolution in the Itô framework

$$dy_t = f(y_t)\,dt + \sigma\, g(y_t)\,dW_t.$$

where $y_t \in U \subset R,\ f: U \to R,\ g: U \to R,$ and $W_t$ is a scalar Brownian motion process. Suppose the Carleman linearization order is three. Then the 'associated' Carleman linearized Markovian stochastic evolution in the Itô framework under the 'finite closure' is

$$d\begin{pmatrix} y_t \\ y_t^2 \\ y_t^3 \end{pmatrix} = \left( \begin{pmatrix} f \\ \sigma^2 g^2 \\ 0 \end{pmatrix} + \begin{pmatrix} f' & \frac{1}{2!} f'' & \frac{1}{3!} f''' \\ 2f + 2\sigma^2 g\, g' & 2f' + \sigma g'^2 + 2\sigma^2 g\, g'' & f'' + \sigma^2 g' g'' + \frac{2}{3!}\sigma^2 g\, g''' \\ 3\sigma^2 g^2 & 3f + 6\sigma^2 g\, g' & 3f' + 3\sigma g'^2 + 6\sigma^2 g\, g'' \end{pmatrix} \begin{pmatrix} y_t \\ y_t^2 \\ y_t^3 \end{pmatrix} \right) dt$$

$$+ \begin{pmatrix} \sigma g' & \frac{1}{2!}\sigma g'' & \frac{1}{3!}\sigma g''' \\ 2\sigma g & 2\sigma g' & \sigma g'' \\ 0 & 3\sigma g & 3\sigma g' \end{pmatrix} \begin{pmatrix} y_t \\ y_t^2 \\ y_t^3 \end{pmatrix} dW_t + \begin{pmatrix} \sigma g \\ 0 \\ 0 \end{pmatrix} dW_t,$$

where $f = f(y_t)$, $f' = \frac{df(y_t)}{dy_t}$, $f'' = \frac{d^2(y_t)}{dy_t^2}$, $f''' = \frac{d^3(y_t)}{dy_t^3}$, $g = g(y_t)$, $g' = \frac{dg(y_t)}{dy_t}$, $g'' = \frac{d^2 g(y_t)}{dy_t^2}$ and $g''' = \frac{d^3 g(y_t)}{dy_t^3}$ evaluated at $y_t = 0$.

*Proof:* Here we weave the proof of the Theorem via illustrating the Carleman embedding into the above SDE of the Theorem. Using the expansions of the functions $f(y_t)$ and $g(y_t)$ from the generating function perspective and the function analyticity, we have

$$dy_t = \sum_{0 \le k} a_k y_t^k dt + \sum_{0 \le k} b_k y_t^k dW_t,$$

where

$$a_k = \frac{1}{k!}\frac{d^k f(y_t)}{dy_t^k}|y_t = 0, \;\; b_k = \frac{\sigma}{k!}\frac{d^k g(y_t)}{dy_t^k}|y_t = 0.$$

Suppose the contribution beyond the cubic nonlinearity is smaller, thus

$$dy_t \approx (a_0 + a_1 y_t + a_2 y_t^2 + a_3 y_t^3)dt + (b_0 + b_1 y_t + b_2 y_t^2 + b_3 y_t^3)dW_t.$$

Thanks to the stochastic differential rule (Karatzas and Shreve [7], p. 148), we have the stochastic evolution

$$\begin{aligned} dy_t^2 &= 2y_t dy_t + (dy_t)^2 \\ &= 2y_t((a_0 + a_1 y_t + a_2 y_t^2 + a_3 y_t^3)dt + (b_0 + b_1 y_t + b_2 y_t^2 + b_3 y_t^3)dW_t) \\ &\quad + (b_0 + b_1 y_t + b_2 y_t^2 + b_3 y_t^3)^2 dt \\ &= (b_0^2 + (2a_0 + 2b_0 b_1)y_t + (2a_1 + b_1^2 + 2b_0 b_2)y_t^2 + (2a_2 + 2b_1 b_2 + 2b_0 b_3)y_t^3 \\ &\quad + (2a_3 + b_2^2 + 2b_1 b_3)y_t^4 + 2b_2 b_3 y_t^5 + b_3^2 y_t^6)dt + (2b_0 y_t + 2b_1 y_t^2 + 2b_2 y_t^3 + 2b_3 y_t^4)dW_t. \end{aligned}$$

$$\begin{aligned} dy_t^3 &= y_t\, dy_t^2 + y_t^2\, dy_t + dy_t\, dy_t^2 \\ &= y_t((b_0^2 + (2a_0 + 2b_0 b_1)y_t + (2a_1 + b_1^2 + 2b_0 b_2)y_t^2 + (2a_2 + 2b_1 b_2 + 2b_0 b_3)y_t^3 \\ &\quad + (2a_3 + b_2^2 + 2b_1 b_3)y_t^4 + 2b_2 b_3 y_t^5 + b_3^2 y_t^6)dt + (2b_0 y_t + 2b_1 y_t^2 + 2b_2 y_t^3 + 2b_3 y_t^4)dW_t) \\ &\quad + y_t^2(a_0 + a_1 y_t + a_2 y_t^2 + a_3 y_t^3)dt + (b_0 y_t^2 + b_1 y_t^3 + b_2 y_t^4 + b_3 y_t^5)dW_t) \end{aligned}$$

$$+(b_0+b_1y_t+b_2y_t^2+b_3y_t^3)(2b_0y_t+2b_1y_t^2+2b_2y_t^3+2b_3y_t^4)\,dt$$

$$=(3b_0^2y_t+(3a_0+6b_0b_1)y_t^2+(3a_1+3b_1^2+6b_0b_2)y_t^3+(3a_2+6b_1b_2+6b_0b_3)y_t^4$$
$$+(3a_3+3b_2^2+6b_1b_3)y_t^5+6b_2b_3y_t^6+3b_3^2y_t^7)dt+(3b_0y_t^2+3b_1y_t^3+3b_2y_t^4+3b_3y_t^5)dW_t.$$

After adopting the finite closure, i.e. truncating up to the cubic nonlinearity terms, as a consequence of invoking a condition, i.e. the smaller contribution of nonlinearities beyond the cubic nonlinearity to the stochastic evolution (Rugh [8], p. 108), we get

$$dy_t^2\approx(b_0^2+2(a_0+b_0b_1)y_t+(b_1^2+2a_1+2b_0b_2)y_t^2+(2a_2+2b_1b_2+2b_0b_3)y_t^3)dt$$
$$+(2b_0y_t+2b_1y_t^2+2b_2y_t^3)dW_t, \tag{2a}$$

$$dy_t^3\approx(3b_0^2y_t+(3a_0+6b_0b_1)y_t^2+(3a_1+3b_1^2+6b_0b_2)y_t^3)dt+(3b_0y_t^2+3b_1y_t^3)dW_t. \tag{2b}$$

After combining equations (1) and (2a-2b), we have

$$d\begin{pmatrix}y_t\\y_t^2\\y_t^3\end{pmatrix}=\left(\begin{pmatrix}a_0\\b_0^2\\0\end{pmatrix}+\begin{pmatrix}a_1&a_2&a_3\\2a_0+2b_0b_1&2a_1+b_1^2+2b_0b_2&2a_2+2b_1b_2+2b_0b_3\\3b_0^2&3a_0+6b_0b_1&3a_1+3b_1^2+6b_0b_2\end{pmatrix}\begin{pmatrix}y_t\\y_t^2\\y_t^3\end{pmatrix}\right)dt$$
$$+\begin{pmatrix}b_1&b_2&b_3\\2b_0&2b_1&2b_2\\0&3b_0&3b_1\end{pmatrix}\begin{pmatrix}y_t\\y_t^2\\y_t^3\end{pmatrix}dW_t+\begin{pmatrix}b_0\\0\\0\end{pmatrix}dW_t. \tag{3}$$

An alternative notational rephrasing is

$$d\begin{pmatrix}y_1(t)\\y_2(t)\\y_3(t)\end{pmatrix}=\left(\begin{pmatrix}a_0\\b_0^2\\0\end{pmatrix}+\begin{pmatrix}a_1&a_2&a_3\\2a_0+2b_0b_1&2a_1+b_1^2+2b_0b_2&2a_2+2b_1b_2+2b_0b_3\\3b_0^2&3a_0+6b_0b_1&3a_1+3b_1^2+6b_0b_2\end{pmatrix}\begin{pmatrix}y_1(t)\\y_2(t)\\y_3(t)\end{pmatrix}\right)dt$$
$$+\begin{pmatrix}b_1&b_2&b_3\\2b_0&2b_1&2b_2\\0&3b_0&3b_1\end{pmatrix}\begin{pmatrix}y_1(t)\\y_2(t)\\y_3(t)\end{pmatrix}dW_t+\begin{pmatrix}b_0\\0\\0\end{pmatrix}dW_t,$$

where $(y_t\quad y_t^2\quad y_t^3)^T=(y_1(t)\quad y_2(t)\quad y_3(t))^T$. The above can be further recast as

$$d\xi_t=(A_0+A_t\xi_t)dt+(D_t\xi_t+L_t)dW_t,$$

where

$$\xi_t=\begin{pmatrix}y_t\\y_t^2\\y_t^3\end{pmatrix},\ A_0=\begin{pmatrix}a_0\\b_0^2\\0\end{pmatrix},\quad A_t=\begin{pmatrix}a_1&a_2&a_3\\2a_0+2b_0b_1&2a_1+b_1^2+2b_0b_2&2a_2+2b_1b_2+2b_0b_3\\3b_0^2&3a_0+6b_0b_1&3a_1+3b_1^2+6b_0b_2\end{pmatrix},$$

$$D_t=\begin{pmatrix}b_1&b_2&b_3\\2b_0&2b_1&2b_2\\0&3b_0&3b_1\end{pmatrix},\ L_t=\begin{pmatrix}b_0\\0\\0\end{pmatrix}.$$

In the direct form, the SDE boils down to

$$d\xi_t = (A_0 + A_t\xi_t)dt + (D_t\xi_t + L_t)dW_t, \tag{4}$$

where

$$\xi_t = \begin{pmatrix} y_t \\ y_t^2 \\ y_t^3 \end{pmatrix},\ A_0 = \begin{pmatrix} f \\ \sigma^2 g^2 \\ 0 \end{pmatrix},\ D_t = \begin{pmatrix} \sigma g' & \frac{1}{2!}\sigma g'' & \frac{1}{3!}\sigma g''' \\ 2\sigma g & 2\sigma g' & \sigma g'' \\ 0 & 3\sigma g & 3\sigma g' \end{pmatrix},\ L_t = \begin{pmatrix} \sigma g \\ 0 \\ 0 \end{pmatrix},$$

$$A_t = \begin{pmatrix} f' & \frac{1}{2!} f'' & \frac{1}{3!} f''' \\ 2f + 2\sigma^2 g g' & 2f' + \sigma g'^2 + 2\sigma^2 g g'' & f'' + \sigma^2 g' g'' + \frac{2}{3!}\sigma^2 g g''' \\ 3\sigma^2 g^2 & 3f + 6\sigma^2 g g' & 3f' + 3\sigma g'^2 + 6\sigma^2 g g'' \end{pmatrix}.$$

Note that the entries of the associated vectors and matrices are evaluated at $y_t = 0.$

*QED*

Remark 1: Equation (4) is a consequence of the Carleman embedding into the nonlinear SDE $dy_t = f(y_t)\,dt + \sigma g(y_t)\,dW_t$. The Carleman linearized stochastic differential equation assumes the structure of a 'bilinear' stochastic differential equation as well as treats the 'nonlinearity' as a state variable. This implies that the Carleman linearized stochastic differential equation (4) and nonlinear stochastic differential equation (1) account for nonlinearity effects. Similar to an advantage of the finite closure, the Carleman linearization circumvents quite the curse of the formidable complexity attributed to the nonlinearity. To achieve the estimation, we adopt the Carleman linearized formalism in our paper to exploit the usefulness of bilinear stochastic differential equations.

# 3. Embedding the Carleman linearization into the machine swing dynamics

Here we explain succinctly the nonlinear machine swing dynamics and then account for the associated stochasticity. As a result of this, we weave a nonlinear stochastic differential equation. Notably, extracting states from noisy observables becomes a filtering problem. In standard literature, first we achieve the nonlinear filtering, then linearization and finally, discretization. In contrast to available methods for filtering, first linearize the stochastic dynamics, then achieve the filtering and finally the discretization. This method of state extractions from noisy observables has found attention in the literature (Bhatt and Sharma [2], Germani *et al.* [9]). Thus, first we embed the Carleman linearization into the machine swing SDE.

In a deterministic setting, the second-order nonlinear differential equation for relative motion of swing equation of a machine swing dynamics (Chen et al. [10]) given as follows

$$M\frac{d^2\delta_t}{dt^2} + D\frac{d\delta_t}{dt} + P_e = P_m. \tag{5}$$

The notations associated with the equation agree with the notations of power system dynamics literature (Lu *et al.* [3], p. 166). For the multi-machine dynamics transient stability can be found in an older, but appealing paper of Willems and Willems [11]. For the brevity of discussions, we omit the detail. The illustrative diagram is given in Fig. 1.

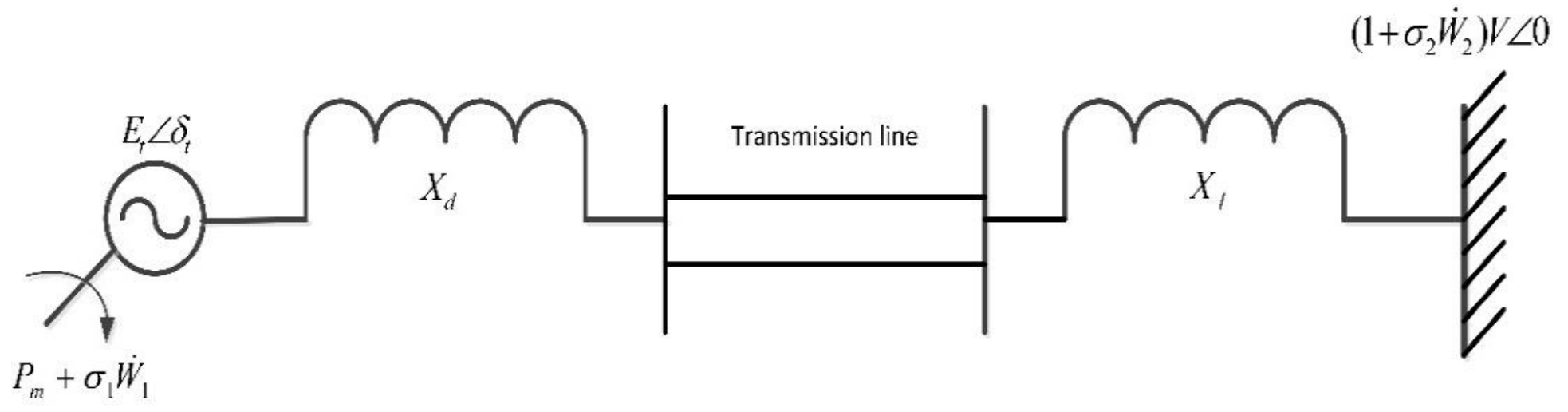

**Fig. 1:** A machine swing dynamic system

Now, we rephrase the swing equation that accounts for noise in the input mechanical power with the additive noise $P_m + \sigma_1 \dot{W}_1(t)$ and the infinite bus voltage magnitude $V(1+\sigma_2 \dot{W}_2(t))$. Thus, we get

$$M\frac{d^2\delta_t}{dt^2} + D\frac{d\delta_t}{dt} + \frac{(1+\sigma_2\dot{W}_2(t))VE_t}{X}\sin\delta_t = P_m + \sigma_1\dot{W}_1(t) \tag{6}$$

The terms $W_1(t)$ and $W_2(t)$ are two independent Wiener processes. The parameters $\sigma_1, \sigma_2$ are the noise intensities of the white noise processes. In phase space formulations, we have

$$dy_t = f(t, y_t)dt + G(t, y_t)dW_t, \tag{7}$$

where

$$y_t = (y_1(t) \quad y_2(t))^T = (\delta_t \quad \dot{\delta}_t)^T \text{ and } f(t, y_t) = (f_1(t, y_1, y_2) \quad f_2(t, y_1, y_2))^T,$$

$$W_t = (W_1(t) \quad W_2(t))^T.$$

Moreover,

$$f(t, y_t) = \begin{pmatrix} y_2 \\ -\frac{D}{M}y_2 - \frac{VE_t}{MX}\sin y_1 + \frac{P_m}{M} \end{pmatrix}, G(t, y_t) = \begin{pmatrix} 0 & 0 \\ \frac{\sigma_1}{M} & -\frac{\sigma_2 VE_t}{MX}\sin y_1 \end{pmatrix}.$$

Here, the term $f(t, y_t)$ shows the nonlinearity in the system and the term $G(t, y_t)$ is the process noise coefficient matrix. The terms $D, M, \sigma_1, \sigma_2, X, E_t, V, P_m$ are the system parameters of the stochastic machine swing dynamic system.

In dynamical systems, power and energy are systematic properties (Willems [12-13]). Consider the power associated with the machine swing is observable in the noisy swing data recording. Thus,

$$dz_t = \left(\frac{VE_t}{X}\right)\sin y_1 dt + d\eta_t, \tag{8}$$

where the $\text{var}(d\eta_t) = \varphi_n dt$ and the state variable $y_1$ represents the rotor angle $\delta_t$. Now next step is to achieve the Carleman embedding into equations (6) and (8).

## 3.1. Carleman linearization of a machine swing SDE

Originally, T. Carleman pioneered a rigorous mathematical technique to transform the sets of polynomial ordinary differential equations into infinite-dimensional linear system representations (Carleman [14], Kawalski and Steeb [15]). Interestingly, we wish to extend the idea of the Carleman linearization for the Markov process and to illustrate the idea, the swing dynamics is analyzed. Its usefulness for the extraction of states from the power system recording data is also unfolded. Here, we take a little pause and explain Carleman linearization for the stochastic process. A connecting thread between the nonlinear stochastic dynamical system and the Carleman linearized stochastic dynamical system. First, we employ the Carleman linearization of the nonlinear SDE. The idea of the Carleman linearization stems from the fact that it decomposes the system nonlinearity and the process noise coefficient in terms of the explicit nonlinearity, e.g. cubic nonlinearity. Secondly, the nonlinearity is treated as a state variable. As a result of these two, we arrive at the Carleman linearized SDE, i.e. bilinear. That is illustrated in Fig. 2, i.e.

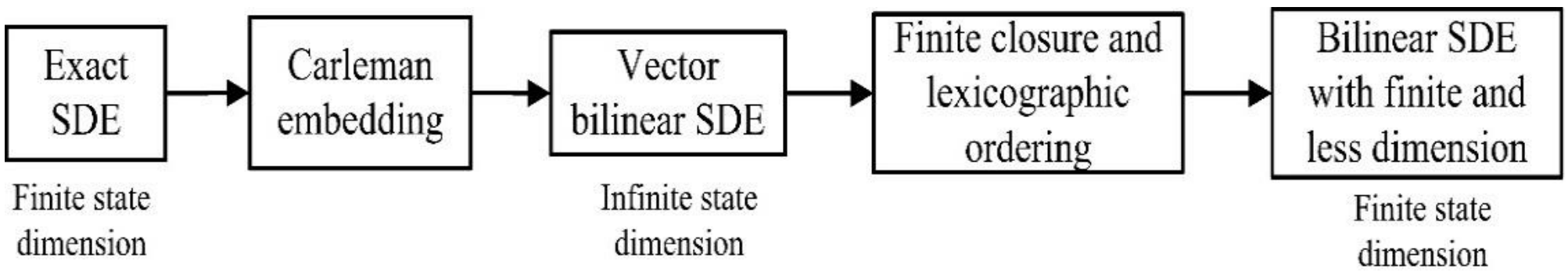

**Fig. 2:** The proposed procedure for embedding the Carleman linearization

The advantage of the Carleman linearization technique is that it accounts for the nonlinearity and embeds bilinearization into the system dynamics. As a result of this, it refines algorithmic procedures as well as simplifies the realization of the concerning algorithmic procedures. For these reasons, the Carleman linearization is the subject of investigation in this paper.

Since the machine swing equation accounts for the sinusoidal nonlinearity, we weave the Carleman embedded in the machine swing equation. Under the finite closure, the Carleman embedding into the stochastic differential equation introduces finite dimensionality into the augmented state vector. The finite closure implies the finite order of the Carleman linearization. The cubic nonlinearity is a generalized nonlinearity (Jing and Lang [6]). Under the finite closure, choose the Carleman linearization order three to account for the cubic nonlinearity. As a result of this, the associated augmented state vector

$$\xi_t = (y_t{}^T \quad y_t^{(2)T} \quad y_t^{(3)T})^T.$$

Suppose

$$y_t = (y_1 \quad y_2)^T, \; y_t^{(2)} = y_t \otimes y_t, \; y_t^{(3)} = y_t \otimes y_t \otimes y_t.$$

The symbol $\otimes$ in the above state represents the Kronecker product. The lexicographic ordering circumvents the problem of redundancy. Thus, in the lexicographic ordering, the components of the augmented state vector are

$$y_t = (y_1 \quad y_2)^T, \; y_t^{(2)} = (y_1^2 \quad y_1 y_2 \quad y_2^2)^T, \; y_t^{(3)} = (y_1^3 \quad y_1^2 y_2 \quad y_1 y_2^2 \quad y_2^3)^T,$$

$$\xi_t = (y_1 \quad y_2 \quad y_1^2 \quad y_1 y_2 \quad y_2^2 \quad y_1^3 \quad y_1^2 y_2 \quad y_1 y_2^2 \quad y_2^3)^T.$$

In a more general setting, suppose the dimension of the state vector $y_t$ is $n$ and the Kronecker power $r$. As a consequence of the lexicographic ordering, the dimension of the vector $y_t^{(r)}$ is $\binom{n+r-1}{r}$. Thanks to the surprising power of the Itô calculus, we sketch the stochastic evolutions of the components of the augmented state vector $\xi_t$ of the paper. Thus, by applying the Itô stochastic differential rule (Pugachev and Sinitsyn [16], Karatzas and Shreve [7]) for the vector $\xi_t$, we obtain the following evolution equations:

$$dy_1 = y_2 dt, \tag{9a}$$

$$dy_2 = \left(-\frac{VE_t}{MX} y_1 - \frac{D}{M} y_2 + \frac{VE_t}{MX} y_1^3 + \frac{P_m}{M}\right) dt + \frac{\sigma_1}{M} dW_1(t)$$
$$-\left(\frac{\sigma_2 VE_t}{MX} y_1 - \frac{\sigma_2 VE_t}{MX} y_1^3\right) dW_2(t), \tag{9b}$$

The augmented evolutions are

$$dy_1^2 = 2 y_1 y_2 dt, \tag{10a}$$

$$dy_1y_2 = \left(\frac{P_m}{M}y_1 - \frac{VE_t}{MX}y_1^2 - \frac{D}{M}y_1y_2 + y_2^2\right)dt + \frac{\sigma_1}{M}y_1dW_1(t)$$

$$-\frac{\sigma_2 VE_t}{MX}y_1^2 dW_2(t), \tag{10b}$$

$$dy_2^2 = \left(\frac{2P_m}{M}y_2 - \frac{2VE_t}{MX}y_1y_2 - \frac{2D}{M}y_2^2 + \frac{\sigma_1^2}{M^2} + \frac{\sigma_2^2V^2E_t^2}{M^2X^2}\right)dt + \frac{2\sigma_1}{M}y_2dW_1(t)$$

$$-\frac{2\sigma_2 VE_t}{MX}y_1y_2dW_2(t), \tag{10c}$$

$$dy_1^3 = 3y_1^2y_2dt, \tag{10d}$$

$$dy_1^2y_2 = \left(\frac{P_m}{M}y_1^2 - \frac{VE_t}{MX}y_1^3 - \frac{D}{M}y_1^2y_2 + 2y_1y_2^2\right)dt + \frac{\sigma_1}{M}y_1^2dW_1(t)$$

$$-\frac{\sigma_2 VE_t}{MX}y_1^3dW_2(t), \tag{10e}$$

$$dy_1y_2^2 = \left(\frac{\sigma_1^2}{M^2}y_1 + \frac{2P_m}{M}y_1y_2 + \frac{\sigma_2^2V^2E_t^2}{M^2X^2}y_1^3 - \frac{2VE_t}{MX}y_1^2y_2 - \frac{2D}{M}y_1y_2^2 + y_2^3\right)dt$$

$$+\frac{2\sigma_1}{M}y_1y_2dW_1(t) - \frac{2\sigma_2VE_t}{MX}y_1^2y_2dW_2(t), \tag{10f}$$

$$dy_2^3 = \left(-3\frac{\sigma_1^2}{M^2}y_2 + \frac{3P_m}{M}y_2^2 - 3\frac{\sigma_2^2V^2E_t^2}{M^2X^2}y_1^2y_2 - \frac{3VE_t}{MX}y_1y_2^2 - \frac{3D}{M}y_2^3\right)dt$$

$$+\frac{3\sigma_1}{M}y_2^2dW_1(t) - \frac{3\sigma_2VE_t}{MX}y_1y_2^2dW_2(t).$$

(10g).

After a rearrangement of the set of equations (9)-(10) in the SDE framework, interestingly, the above set boils down to the following SDE set up:

$$d\xi_t = (A_0 + A_t\xi_t)dt + \sum_{1\le\phi\le2} D_\phi\xi_t dW_\phi + \sum_{1\le\phi\le2} L_\phi dW_\phi, \tag{11}$$

Here, the state vector $\xi_t = (\xi_i) = (y_t^{(r)})$, the notation $r$ denotes the Kronecker power and $1 \le r \le N,\ \ 1 \le i \le \sum_{1\le r\le N}\binom{n+r-1}{r}$, where, $\xi_t \in R^{\sum_{1\le r\le N}\binom{n+r-1}{r}}$. More specifically, the machine swing SDE in the Carleman setting, the term $n = 2,\ \ N = 3, 1 \le r \le 3.$ The partitioned vectors and partitioned matrices, associated with the machine swing SDE in the Carleman setting, i.e. equation (11), are

$$A_0(t) = \begin{pmatrix} A_{01}(t) \\ A_{02}(t) \\ A_{03}(t) \end{pmatrix},\ A_t = \begin{pmatrix} A_{11}(t) & A_{12}(t) & A_{13}(t) \\ A_{21}(t) & A_{22}(t) & A_{23}(t) \\ A_{31}(t) & A_{32}(t) & A_{33}(t) \end{pmatrix},\ D_{1t} = \begin{pmatrix} D_{1t}(1,1) & D_{1t}(1,2) & D_{1t}(1,3) \\ D_{1t}(2,1) & D_{1t}(2,2) & D_{1t}(2,3) \\ D_{1t}(3,1) & D_{1t}(3,2) & D_{1t}(3,3) \end{pmatrix},$$

$$D_{2t}=\begin{pmatrix} D_{2t}(1,1) & D_{2t}(1,2) & D_{2t}(1,3) \\ D_{2t}(2,1) & D_{2t}(2,2) & D_{2t}(2,3) \\ D_{2t}(3,1) & D_{2t}(3,2) & D_{2t}(3,3) \end{pmatrix},\ L_{1t}=\begin{pmatrix} L_{1t}(1) \\ L_{1t}(2) \\ L_{1t}(3) \end{pmatrix},\ L_{2t}=\begin{pmatrix} L_{2t}(1) \\ L_{2t}(2) \\ L_{2t}(3) \end{pmatrix},$$

where

$$A_{01}(t)=\begin{pmatrix} 0 & \frac{P_m}{M} \end{pmatrix}^T,\ A_{02}(t)=\begin{pmatrix} 0 & 0 & \frac{\sigma_1^2}{M^2} \end{pmatrix}^T,\ A_{11}(t)=\begin{pmatrix} 0 & 1 \\ -\frac{VE_t}{MX} & -\frac{D}{M} \end{pmatrix},$$

$$A_{13}(t)=\begin{pmatrix} 0 & 0 & 0 & 0 \\ \frac{VE_t}{6MX} & 0 & 0 & 0 \end{pmatrix},\ A_{21}(t)=\begin{pmatrix} 0 & \frac{P_m}{M} & 0 \\ 0 & 0 & \frac{2P_m}{M} \end{pmatrix}^T,$$

$$A_{22}(t)=\begin{pmatrix} 0 & 2 & 0 \\ -\frac{VE_t}{MX} & -\frac{D}{M} & 1 \\ \frac{\sigma_2^2V^2E_t^2}{M^2X^2} & -\frac{2VE_t}{MX} & -\frac{2D}{M} \end{pmatrix},\ A_{31}(t)=\begin{pmatrix} 0 & 0 & \frac{\sigma_1^2}{M^2} & 0 \\ 0 & 0 & 0 & \frac{3\sigma_1^2}{M^2} \end{pmatrix}^T,$$

$$A_{32}(t)=\begin{pmatrix} 0 & \frac{P_m}{M} & 0 & 0 \\ 0 & 0 & \frac{2P_m}{M} & 0 \\ 0 & 0 & 0 & \frac{3P_m}{M} \end{pmatrix}^T,\ A_{33}(t)=\begin{pmatrix} 0 & 3 & 0 & 0 \\ -\frac{VE_t}{MX} & -\frac{D}{M} & 2 & 0 \\ \frac{\sigma_2^2V^2E_t^2}{M^2X^2} & -\frac{2VE_t}{MX} & -\frac{2D}{M} & 1 \\ 0 & \frac{3\sigma_2^2V^2E_t^2}{M^2X^2} & -\frac{3VE_t}{MX} & -\frac{3D}{M} \end{pmatrix},$$

$$D_{1t}(2,1)=\begin{pmatrix} 0 & \frac{\sigma_1}{M} & 0 \\ 0 & 0 & \frac{2\sigma_1}{M} \end{pmatrix}^T,\ D_{1t}(3,2)=\begin{pmatrix} 0 & \frac{\sigma_1}{M} & 0 & 0 \\ 0 & 0 & \frac{2\sigma_1}{M} & 0 \\ 0 & 0 & 0 & \frac{3\sigma_1}{M} \end{pmatrix}^T,$$

$$D_{2t}(1,1)=\begin{pmatrix} 0 & 0 \\ -\frac{\sigma_2VE_t}{MX} & 0 \end{pmatrix},\ D_{2t}(1,3)=\begin{pmatrix} 0 & 0 & 0 & 0 \\ \frac{\sigma_2VE_t}{6MX} & 0 & 0 & 0 \end{pmatrix},\ L_{1t}(1)=\begin{pmatrix} 0 \\ \frac{\sigma_1}{M} \end{pmatrix},$$

$$D_{2t}(2,2)=\begin{pmatrix} 0 & 0 & 0 \\ -\frac{\sigma_2VE_t}{MX} & 0 & 0 \\ 0 & -\frac{2\sigma_2VE_t}{MX} & 0 \end{pmatrix},\ D_{2t}(3,3)=\begin{pmatrix} 0 & 0 & 0 & 0 \\ -\frac{\sigma_2VE_t}{MX} & 0 & 0 & 0 \\ 0 & \frac{2\sigma_2VE_t}{MX} & 0 & 0 \\ 0 & 0 & \frac{3\sigma_2VE_t}{MX} & 0 \end{pmatrix}.$$

The associated zero matrices of the appropriate sizes are $A_{03}(t)$, $A_{12}(t)$, $A_{23}(t)$, $D_{1t}(1,1)$, $D_{1t}(1,2)$, $D_{1t}(1,3)$, $D_{1t}(2,2)$, $D_{1t}(2,3)$, $D_{1t}(3,1)$, $D_{1t}(3,3)$, $D_{2t}(1,2)$, $D_{2t}(2,1)$, $D_{2t}(2,3)$, $D_{2t}(3,1)$, $D_{2t}(3,2)$, $L_{1t}(2)$, $L_{1t}(3)$, $L_{2t}(1)$, $L_{2t}(2)$ and $L_{2t}(3)$.

The Carleman embedding into the nonlinear noisy observation equation results in the linear framework. Thus,

$$dz_t = C_t \xi_t dt + d\eta_t, \tag{13}$$

where

$$C_t = (C_{1t} \quad C_{2t} \quad C_{3t}),\ C_{1t} = \begin{pmatrix} \frac{VE_t}{X} & 0 \end{pmatrix},\ C_{2t} = \begin{pmatrix} 0 & 0 & 0 \end{pmatrix},\ C_{3t} = \begin{pmatrix} -\frac{VE_t}{6X} & 0 & 0 & 0 \end{pmatrix}.$$

Equation (13) in conjunction with equation (11) is a convenient setup to achieve 'filtering in the Carleman setting'. First, we do filtering in the Carleman setting for extracting states from the machine swing noisy data, then we do in the EKF setting.

## 3.2 Carleman linearized filtering of a machine swing SDE

The filtering of the Carleman linearized stochastic differential system exploits the conditional expectation that accounts for noisy observables accumulated up to the instant $t$ and results in the most probable trajectory of the dynamic system at the instant $t$. Consider a Carleman linearized Itô bilinear stochastic differential equation associated with the linear observation equation, e.g. equations (11)-(13). Then, the stochastic differential equations for filtering in the Carleman setting are

$$d\hat{\xi}_t = (A_0(t) + A_t\hat{\xi}_t)dt + P_t C_t^T \varphi_n^{-1}(dz_t - C_t\hat{\xi}_t dt), \tag{14a}$$

$$dP_t = (P_t A_t^T + A_t P_t + L_t L_t^T + L_t \hat{\xi}_t^T D_t^T + D_t \hat{\xi}_t L_t^T + D_t P_t D_t^T$$
$$+ D_t\, \hat{\xi}_t \hat{\xi}_t^T D_t^T - P_t C_t^T \varphi_n^{-1} C_t P_t)dt, \tag{14b}$$

where

$$\hat{\xi}_t = \begin{pmatrix} \hat{y}_t \\ \hat{y}_t^{(2)} \\ \vdots \\ \hat{y}_t^{(N)} \end{pmatrix},\quad P_t = E(\xi_t - \hat{\xi}_t)(\xi_t - \hat{\xi}_t)^T = \begin{pmatrix} P_{y_t y_t} & P_{y_t y_t^{(2)}} & \cdots & P_{y_t y_t^{(N)}} \\ P_{y_t^{(2)} y_t} & P_{y_t^{(2)} y_t} & \cdots & P_{y_t^{(2)} y_t^{(N)}} \\ \vdots & \vdots & \ddots & \vdots \\ P_{y_t^{(N)} y_t} & P_{y_t^{(N)} y_t^{(2)}} & \cdots & P_{y_t^{(N)} y_t^{(N)}} \end{pmatrix}.$$

Note that the expectation operator $E$ is conditional on accumulated noisy observables. The higher the order of Carleman linearization would be the greater the size of the associated partitioned matrices. As a specific case for practical problems, consider the Carleman order three, then equations (14a)-(14b) have the following stochastic evolutions:

$$d\hat{y}_t = (A_{01}(t) + A_{11}(t)\hat{y}_t + A_{12}(t)\hat{y}_t^{(2)} + A_{13}(t)\hat{y}_t^{(3)})dt$$
$$+ (P_{y_t y_t} C_{1t}^T + P_{y_t y_t^{(2)}} C_{2t}^T + P_{y_t y_t^{(3)}} C_{3t}^T)\varphi_n^{-1}(dz_t - (C_{1t}\hat{y}_t + C_{2t}\hat{y}_t^{(2)} + C_{3t}\hat{y}_t^{(3)})dt), \tag{15a}$$

$$d\hat{y}_t^{(2)} = (A_{02}(t) + A_{21}(t)\hat{y}_t + A_{22}(t)\hat{y}_t^{(2)} + A_{23}(t)\hat{y}_t^{(3)})dt$$
$$+ (P_{y_t^{(2)} y_t} C_{1t}^T + P_{y_t^{(2)} y_t^{(2)}} C_{2t}^T + P_{y_t^{(2)} y_t^{(3)}} C_{3t}^T)\varphi_n^{-1}(dz_t - (C_{1t}\hat{y}_t + C_{2t}\hat{y}_t^{(2)} + C_{3t}\hat{y}_t^{(3)})dt), \tag{15b}$$

$$d\hat{y}_t^{(3)} = (A_{03}(t) + A_{31}(t)\hat{y}_t + A_{32}(t)\hat{y}_t^{(2)} + A_{33}(t)\hat{y}_t^{(3)})dt$$
$$+ (P_{y_t^{(3)} y_t} C_{1t}^T + P_{y_t^{(3)} y_t^{(2)}} C_{2t}^T + P_{y_t^{(3)} y_t^{(3)}} C_{3t}^T)\varphi_n^{-1}(dz_t - (C_{1t}\hat{y}_t + C_{2t}\hat{y}_t^{(2)} + C_{3t}\hat{y}_t^{(3)})dt). \tag{15c}$$

Suppose the dimension of the state vector is two. Then, equations (15a)-(15c) concerning the machine swing SDE (11) and the observation equation (13) assumes the following structures:

$$d\hat{y}_1 = \hat{y}_2\, dt + \left( \frac{VE_t}{X} P_{y_1} - \frac{VE_t}{2X} \hat{y}_1{}^2 P_{y_1} - \frac{VE_t}{2X} P_{y_1}^2 \right) \varphi_n^{-1} \left( dz_t - \left( \frac{VE_t}{X} \hat{y}_1 - \frac{VE_t}{6X} \hat{y}_1{}^3 \right) dt \right), \tag{16a}$$

$$d\hat{y}_2 = \left( \frac{P_m}{M} - \frac{VE_t}{MX} \hat{y}_1 - \frac{D}{M} \hat{y}_2 + \frac{VE_t}{2MX} \hat{y}_1 P_{y_1} + \frac{VE_t}{6MX} \hat{y}_1{}^3 \right) dt$$

$$+ \left( \frac{VE_t}{X} P_{y_1 y_2} - \frac{VE_t}{2X} \hat{y}_1{}^2 P_{y_1 y_2} - \frac{VE_t}{2X} P_{y_1} P_{y_1 y_2} \right) \varphi_n^{-1} \left( dz_t - \left( \frac{VE_t}{X} \hat{y}_1 - \frac{VE_t}{6X} \hat{y}_1{}^3 \right) dt \right), \tag{16b}$$

$$d\hat{y}_1^2 = \left( 2P_{y_1 y_2} + 2\hat{y}_1 \hat{y}_2 \right) dt + \left( \frac{2VE_t}{X} \hat{y}_1 P_{y_1} - \frac{VE_t}{X} \hat{y}_1{}^3 P_{y_1} - \frac{5VE_t}{2X} \hat{y}_1 P_{y_1}^2 \right)$$

$$\varphi_n^{-1} \left( dz_t - \left( \frac{VE_t}{X} \hat{y}_1 - \frac{VE_t}{6X} \hat{y}_1{}^3 \right) dt \right), \tag{16c}$$

$$dE(y_1 y_2) = \left( \frac{P_m}{M} \hat{y}_1 - \frac{VE_t}{MX} P_{y_1} - \frac{VE_t}{MX} \hat{y}_1{}^2 - \frac{D}{M} P_{y_1 y_2} - \frac{D}{M} \hat{y}_1 \hat{y}_2 + P_{y_2} + \hat{y}_2{}^2 \right) dt$$

$$+ \left( \frac{VE_t}{X} \hat{y}_2 P_{y_1} + \frac{VE_t}{X} \hat{y}_1 P_{y_1 y_2} - \frac{VE_t}{2X} \hat{y}_1{}^2 \hat{y}_2 P_{y_1} - \frac{VE_t}{2X} \hat{y}_1{}^3 P_{y_1 y_2} - \frac{3VE_t}{2X} \hat{y}_1 P_{y_1} P_{y_1 y_2} - \frac{VE_t}{2X} \hat{y}_2 P_{y_1}^2 \right)$$

$$\varphi_n^{-1} \left( dz_t - \left( \frac{VE_t}{X} \hat{y}_1 - \frac{VE_t}{6X} \hat{y}_1{}^3 \right) dt \right), \tag{16d}$$

$$d\hat{y}_2^2 = \left( \frac{\sigma_1^2}{M^2} + \frac{2P_m}{M} \hat{y}_2 + \frac{\sigma_2^2 V^2 E_2^2}{M^2 X^2} P_{y_1} + \frac{\sigma_2^2 V^2 E_2^2}{M^2 X^2} \hat{y}_1{}^2 - \frac{2VE_t}{MX} P_{y_1 y_2} - \frac{2VE_t}{MX} \hat{y}_1 \hat{y}_2 \right.$$

$$\left. - \frac{2D}{M} P_{y_2} - \frac{2D}{M} \hat{y}_2{}^2 \right) dt + \left( \frac{2VE_t}{X} \hat{y}_2 P_{y_1 y_2} - \frac{VE_t}{X} \hat{y}_1{}^2 \hat{y}_2 P_{y_1 y_2} - \frac{VE_t}{X} \hat{y}_2 P_{y_1} P_{y_1 y_2} \right.$$

$$\left. - \frac{VE_t}{2X} \hat{y}_1 P_{y_1 y_2}^2 \right) \varphi_n^{-1} \left( dz_t - \left( \frac{VE_t}{X} \hat{y}_1 - \frac{VE_t}{6X} \hat{y}_1{}^3 \right) dt \right), \tag{16e}$$

$$d\hat{y}_1^3 = \left( 3\hat{y}_2 P_{y_1} + 3\hat{y}_1{}^2 \hat{y}_2 + 6\hat{y}_1 P_{y_1 y_2} \right) dt + \left( \frac{3VE_t}{X} \hat{y}_1{}^2 P_{y_1} + \frac{3VE_t}{X} P_{y_1}^2 - \frac{3VE}{2X} \hat{y}_1{}^4 P_{y_1} \right.$$

$$\left. - \frac{6VE_t}{X} \hat{y}_1{}^2 P_{y_1}^2 - \frac{5VE_t}{2X} P_{y_1}^3 \right) \varphi_n^{-1} \left( dz_t - \left( \frac{VE_t}{X} \hat{y}_1 - \frac{VE_t}{6X} \hat{y}_1{}^3 \right) dt \right), \tag{16f}$$

$$dE(y_1^2 y_2) = \left( \frac{P_m}{M} P_{y_1} + \frac{P_m}{M} \hat{y}_1{}^2 - \frac{3VE_t}{MX} \hat{y}_1 P_{y_1} - \frac{VE_t}{MX} \hat{y}_1{}^3 - \frac{D}{M} \hat{y}_2 P_{y_1} - \frac{D}{M} \hat{y}_1{}^2 \hat{y}_2 \right.$$

$$\left. - \frac{2D}{M} \hat{y}_1 P_{y_1 y_2} + 2\hat{y}_1 P_{y_2} + 2\hat{y}_1 \hat{y}_2{}^2 + 4\hat{y}_2 P_{y_1 y_2} \right) dt + \left( \frac{2VE_t}{X} \hat{y}_1 \hat{y}_2 P_{y_1} + \frac{VE_t}{X} \hat{y}_1{}^2 P_{y_1 y_2} \right.$$

$$+ \frac{3VE_t}{X} P_{y_1} P_{y_1 y_2} - \frac{VE_t}{X} \hat{y}_1{}^3 \hat{y}_2 P_{y_1} - \frac{VE_t}{2X} \hat{y}_1{}^4 P_{y_1 y_2} - \frac{4VE_t}{X} \hat{y}_1{}^2 P_{y_1} P_{y_1 y_2} \tag{16g}$$

$$-\frac{2VE_t}{X}\hat{y}_1\hat{y}_2P_{y_1}^2-\frac{3VE}{2X}P_{y_1}^2P_{y_1y_2}\Bigg)\varphi_n^{-1}\left(dz_t-\left(\frac{VE_t}{X}\hat{y}_1-\frac{VE_t}{6X}\hat{y}_1^{\,3}\right)dt\right),$$

$$dE(y_1y_2^2)=\Bigg(\frac{\sigma_1^2}{M^2}\hat{y}_1+\frac{2P_m}{M}P_{y_1y_2}+\frac{2P_m}{M}\hat{y}_1\hat{y}_2+\frac{3\sigma_2^2V^2E_t^2}{M^2X^2}\hat{y}_1P_{y_1}+\frac{\sigma_2^2V^2E_t^2}{M^2X^2}\hat{y}_1^{\,3}$$

$$-\frac{2VE_t}{MX}\hat{y}_2P_{y_1}-\frac{2VE_t}{MX}\hat{y}_1^{\,2}\hat{y}_2-\frac{4VE_t}{MX}\hat{y}_1P_{y_1y_2}-\frac{2D}{M}\hat{y}_1P_{y_2}-\frac{2D}{M}\hat{y}_1\hat{y}_2^{\,2}-\frac{4D}{M}\hat{y}_2P_{y_1y_2}+3\hat{y}_2P_{y_2}$$

$$+\hat{y}_2^{\,3})dt+\Bigg(\frac{VE_t}{X}\hat{y}_2^{\,2}P_{y_1}+\frac{2VE_t}{X}\hat{y}_1\hat{y}_2P_{y_1y_2}+\frac{VE_t}{X}P_{y_1}P_{y_2}+\frac{VE_t}{X}P_{y_1y_2}^2$$

$$-\frac{3VE_t}{X}\hat{y}_1\hat{y}_2P_{y_1}P_{y_1y_2}-\frac{VE_t}{2X}\hat{y}_1^{\,2}P_{y_1}P_{y_2}-\frac{VE_t}{X}\hat{y}_1^{\,2}P_{y_1y_2}^2-\frac{VE_t}{2X}\hat{y}_2^{\,2}P_{y_1}^2$$

$$-\frac{VE_t}{2X}P_{y_1}^2P_{y_2}-\frac{2VE_t}{X}P_{y_1}P_{y_1y_2}^2\Bigg)\varphi_n^{-1}\left(dz_t-\left(\frac{VE_t}{X}\hat{y}_1-\frac{VE_t}{6X}\hat{y}_1^{\,3}\right)dt\right),$$

(16h)

$$d\hat{y}_2^3=\Bigg(\frac{3\sigma_1^2}{M^2}\hat{y}_2+\frac{3P_m}{M}P_{y_2}+\frac{3P_m}{M}\hat{y}_2^{\,2}+\frac{3\sigma_2^2V^2E_t^2}{M^2X^2}\hat{y}_2P_{y_1}+\frac{3\sigma_2^2V^2E_t^2}{M^2X^2}\hat{y}_1^{\,2}\hat{y}_2$$

$$+\frac{6\sigma_2^2V^2E_t^2}{M^2X^2}\hat{y}_1P_{y_1y_2}-\frac{3VE_t}{MX}\hat{y}_1P_{y_2}-\frac{3VE_t}{MX}\hat{y}_1\hat{y}_2^{\,2}-\frac{6VE_t}{MX}\hat{y}_2P_{y_1y_2}-\frac{9D}{M}\hat{y}_2P_{y_2}$$

$$-\frac{3D}{M}\hat{y}_2^{\,3}\Bigg)dt+\Bigg(\frac{3VE_t}{X}\hat{y}_1^{\,2}P_{y_1y_2}+\frac{3VE_t}{X}P_{y_2}P_{y_1y_2}-\frac{3VE_t}{2X}P_{y_1}P_{y_1y_2}P_{y_2}-\frac{VE_t}{X}P_{y_1y_2}^3$$

$$-\frac{3VE_t}{2X}\hat{y}_2^{\,2}P_{y_1}P_{y_1y_2}-\frac{3VE_t}{2X}\hat{y}_1\hat{y}_2P_{y_1y_2}^2-\frac{3VE_t}{2X}\hat{y}_1^2P_{y_2}P_{y_1y_2}-\hat{y}_1^2\hat{y}_2^2\frac{3VE_t}{2X}P_{y_1y_2}\Bigg)$$

$$\varphi_n^{-1}\left(dz_t-\left(\frac{VE_t}{X}\hat{y}_1-\frac{VE_t}{6X}\hat{y}_1^{\,3}\right)dt\right).$$

(16i)

The Carleman linearized machine swing conditional variance evolution equations are obtained by using equations (11) and (13) in combination with the filtering equation (14b). As a result of this, we obtain the following stochastic equations:

$$dP_{y_1}=\left(2P_{y_1y_2}+\frac{2VE_t}{MX}\hat{y}_1P_{y_1}^2-\left(\frac{VE_t}{X}P_{y_1}-\frac{VE_t}{2X}\hat{y}_1^{\,2}P_{y_1}+\frac{VE_t}{2X}P_{y_1}^2\right)^2\varphi_n^{-1}\right)dt,$$

(17a)

$$dP_{y_2}=\Bigg(\frac{\sigma_1^2}{M^2}+\frac{\sigma_2^2V^2E_t^2}{M^2X^2}P_{y_1}+\frac{\sigma_2^2V^2E_t^2}{M^2X^2}\hat{y}_1^2-2\frac{VE_t}{MX}P_{y_1y_2}-2\frac{D}{M}P_{y_2}-\frac{VE_t}{MX}\hat{y}_1\hat{y}_2P_{y_1}$$

$$-\frac{VE_t}{3MX}\hat{y}_1^{\,2}\hat{y}_2+\frac{VE_t}{MX}\hat{y}_1P_{y_1}P_{y_2}+\frac{VE_t}{2MX}\hat{y}_1P_{y_1y_2}^2-\left(\frac{VE_t}{X}P_{y_1y_2}-\frac{VE_t}{2X}\hat{y}_1^{\,2}P_{y_1y_2}+\frac{VE_t}{2X}P_{y_1}P_{y_1y_2}\right)^2\varphi_n^{-1}\Bigg)dt,$$

(17b)

$$dP_{y_1y_2}=\Bigg(-\frac{VE_t}{MX}P_{y_1}-\frac{D}{M}P_{y_1y_2}+P_{y_2}+\frac{2VE_t}{MX}\hat{y}_1P_{y_1}P_{y_1y_2}-\frac{VE_t}{2MX}\hat{y}_1^{\,2}P_{y_1}-\frac{VE_t}{6MX}\hat{y}_1^{\,4}$$

(17c)

$$-\left(\frac{VE_t}{X}P_{y_1}-\frac{VE_t}{2X}\hat{y}_1{}^2P_{y_1}+\frac{VE_t}{2X}P_{y_1}^2\right)\left(\frac{VE_t}{X}P_{y_1y_2}-\frac{VE_t}{2X}\hat{y}_1{}^2P_{y_1y_2}+\frac{VE_t}{2X}P_{y_1}P_{y_1y_2}\right)\varphi_n^{-1}\Bigg)dt.$$

Remark 2: the notation $E$ stands for the conditional expectation. For the brevity of notations, $Ey_t^n=\hat{y}_t^n$ and $(Ey_t)^n=\hat{y}_t{}^n$. The relationship between $Ey_t^n$ and $(Ey_t)^n$ under nearly Gaussian assumption can be found in *appendix* of the paper. For the Carleman linearization order two, the machine swing filtering becomes the linear that is of little interest. The machine swing nonlinear filtering equations of the paper adopt the Carleman linearization order three that is of greater interest.

The stochastic nonlinear machine swing equation describes the rotor dynamics of the synchronous machines and it helps to stabilize the power system. It is imperative to develop filtering algorithms for extracting the states from the noisy machine swing data. That would useful to anticipate power system operations as well as gaining an idea about the fragility of the power system dynamics. Since the Carleman linearization-based filtering of the states has superiority in the senses of the better estimates as well as more simplified realization of the filters in contrast to nonlinear filters. Now, we write the proposed algorithmic estimation procedure for filtering of machine swing dynamic in the Carleman setting.

Step 1 Adding suitable noise terms into the nonlinear machine swing dynamics

$$dy_t=f(y_t)dt.$$

Step 2 Reorganizing the nonlinear system dynamics of step 1 that produces in the form of Itô stochastic differential equations of the machine swing dynamics.

Step 3 Accomplishing the power series expansion of the system nonlinearity and the process noise terms for the Carleman linearization order $N,\ N>1$.

Step 4 Reorganizing the equation of step 3 using the partitioning of the associated states and matrices, which results in the Carleman linearized bilinear stochastic differential equation:

$$d\xi_t=(A_0(t)+A_t\xi_t)dt+D_t\xi_t dW_t+L_tW_t.$$

Step 5 Reorganize the associated nonlinear observation equation into the linear observation equation using the Carleman linearization of step similar to step 3, i.e.

$$dz_t=C_t\xi_t dt+d\eta_t.$$

Step 6 Generate the Carleman linearized filtered estimates for the machine swing dynamics by exploiting the filtering equations (16a)-(16i) and (17a)-(17c).

### 3.3 Extended Kalman filtering of a machine swing SDE

In estimation theory, the benchmark EKF approach is the nonlinear form of the Kalman filter, which is linearized about an estimate of the Markovian states. The Kalman filter is the optimal linear estimator for linear system models with additive independent white noise in both the transition and the measurement systems. To contrast the Carleman linearized filtered states with the benchmark method, we exploit the celebrated EKF. Notably, the nonlinear filtering of the machine swing dynamics is not sufficiently well known in the literature.

The general moment equations of the extended Kalman filter for the continuous state-continuous measurement system, the nonlinear system combined with the nonlinear observation equation, are (Jazwinski [17], p. 338)

$$d\hat{y}_t=f_t(t,\hat{y}_t)dt+P_t\left(\frac{\partial h(t,\hat{y}_t)}{\partial\hat{y}_t}\right)^T\varphi_n^{-1}\big(dz_t-h(t,\hat{y}_t)dt\big), \qquad (18a)$$

$$dP_t=\left(P_t\left(\frac{\partial f(t,\hat{y}_t)}{\partial\hat{y}_t}\right)^T+\frac{\partial f(t,\hat{y}_t)}{\partial\hat{y}_t}P_t+(G_tG_t{}^T)-P_t\left(\frac{\partial h(t,\hat{y}_t)}{\partial\hat{y}_t}\right)^T\varphi_n^{-1}\left(\frac{\partial h(t,\hat{y}_t)}{\partial\hat{y}_t}\right)P_t\right)dt. \qquad (18b)$$

The stochastic evolutions of the conditional expectations of the EKF for the machine swing dynamics is

$$d\hat{y}_1 = \hat{y}_2 dt + \frac{VE_t}{X} P_{11} \cos \hat{y}_1 \, \varphi_n^{-1} \left( dz_t - \frac{VE_t}{X} \sin \hat{y}_1 dt \right), \tag{19a}$$

$$d\hat{y}_2 = \left( -\frac{VE_t}{MX} \sin \hat{y}_1 - \frac{D}{M} \hat{y}_2 + \frac{P_m}{M} \right) dt + \frac{VE_t}{X} P_{12} \cos \hat{y}_1 \, \varphi_n^{-1} \left( dz_t - \frac{VE_t}{X} \sin \hat{y}_1 dt \right). \tag{19b}$$

The above coupled equations, (19a)-(19b), are a consequence of combining equations (7)-(8) in the sense of equation (18a).

Furthermore, the concerning conditional variance evolutions for the swing dynamics boil down to

$$dP_{11} = \left( 2P_{12} - \frac{V^2 E_t^2}{X^2} P_{11}^2 \cos^2 \hat{y}_1 \, \varphi_n^{-1} \right) dt, \tag{20a}$$

$$dP_{22} = \left( -\frac{2VE_t}{MX} P_{12} \cos \hat{y}_1 - \frac{2D}{M} P_{22} + \frac{\sigma_1^2}{M^2} + \frac{\sigma_2^2 V^2 E_t^2}{M^2 X^2} \sin^2 \hat{y}_1 - \frac{V^2 E_t^2}{X^2} P_{12}^2 \cos^2 \hat{y}_1 \, \varphi_n^{-1} \right) dt, \tag{20b}$$

$$dP_{12} = \left( P_{22} - \frac{VE_t}{MX} P_{11} \cos \hat{y}_1 - \frac{D}{M} P_{12} - \frac{V^2 E_t^2}{X^2} P_{11} P_{12} \cos^2 \hat{y}_1 \, \varphi_n^{-1} \right) dt. \tag{20c}$$

# 4. Simulation Results

In this section, two sets of simulation parameters associated with the machine swing dynamics are exploited for the numerical study of the noisy machine swing dynamics in the Carleman setting. To examine the usefulness of the proposed filtering algorithm-based estimation approach, i.e. filtering in the Carleman linearization setting, we execute a MATLAB simulation of the stochastic nonlinear stochastic machine swing dynamics. The machine swing filtering equations in the Carleman setting are contrasted with the machine swing filtering in the EKF setting.

## *4.1 First Set*

We demonstrate the usefulness of the Carleman linearization order three graphically. Consider the first set of system parameters mentioned in Table 1 (Grainger and Stevenson [18]). For filtering analysis via simulation, the associated initial conditions are: $y_1(0) = 1,\ y_2(0) = 2,$ $P_{y_1}(0) = P_{y_1 y_2}(0) = 0$ and $P_{y_2}(0) = 2.$ Fig. 3 shows a comparison between Carleman linearized state trajectories with the nonlinear state trajectories, i.e. nonlinear stochastic differential equations.

**Table 1:** The first set of machine swing dynamics parameters

| Parameters | $D$ | $E_t$ | $V$ | $X$ | $\sigma_1$ | $\sigma_2$ | $P_m$ | $M$ | $\varphi_n$ |
|---|---|---|---|---|---|---|---|---|---|
| Values | 0.25 | 1.2 | 1 | 0.25 | 0.08 | 0.06 | 1 | 0.2 | 100 |
| Units | $pu/rad/s$ | $pu$ | $pu$ | $pu$ | - | - | $pu$ | - | - |

Fig. 3 has two parts, i.e. Fig. 3(a) and Fig. 3(b) for the states $y_1$ and $y_2$ trajectories respectively. The solid lines illustrate the state evolutions, a consequence of the exact machine swing SDE. The marked dash-dash lines show the state evolutions, a consequence of the SDE with the Carleman linearization. A step-change in the rotor angle $\delta_t$ is introduced at the time interval of 5 $\sec,$ i.e.

the absolute disturbance of 2 $rad$. That unfolds the Carleman linearized SDE trajectory tracks the nonlinear SDE trajectory quite very closely.

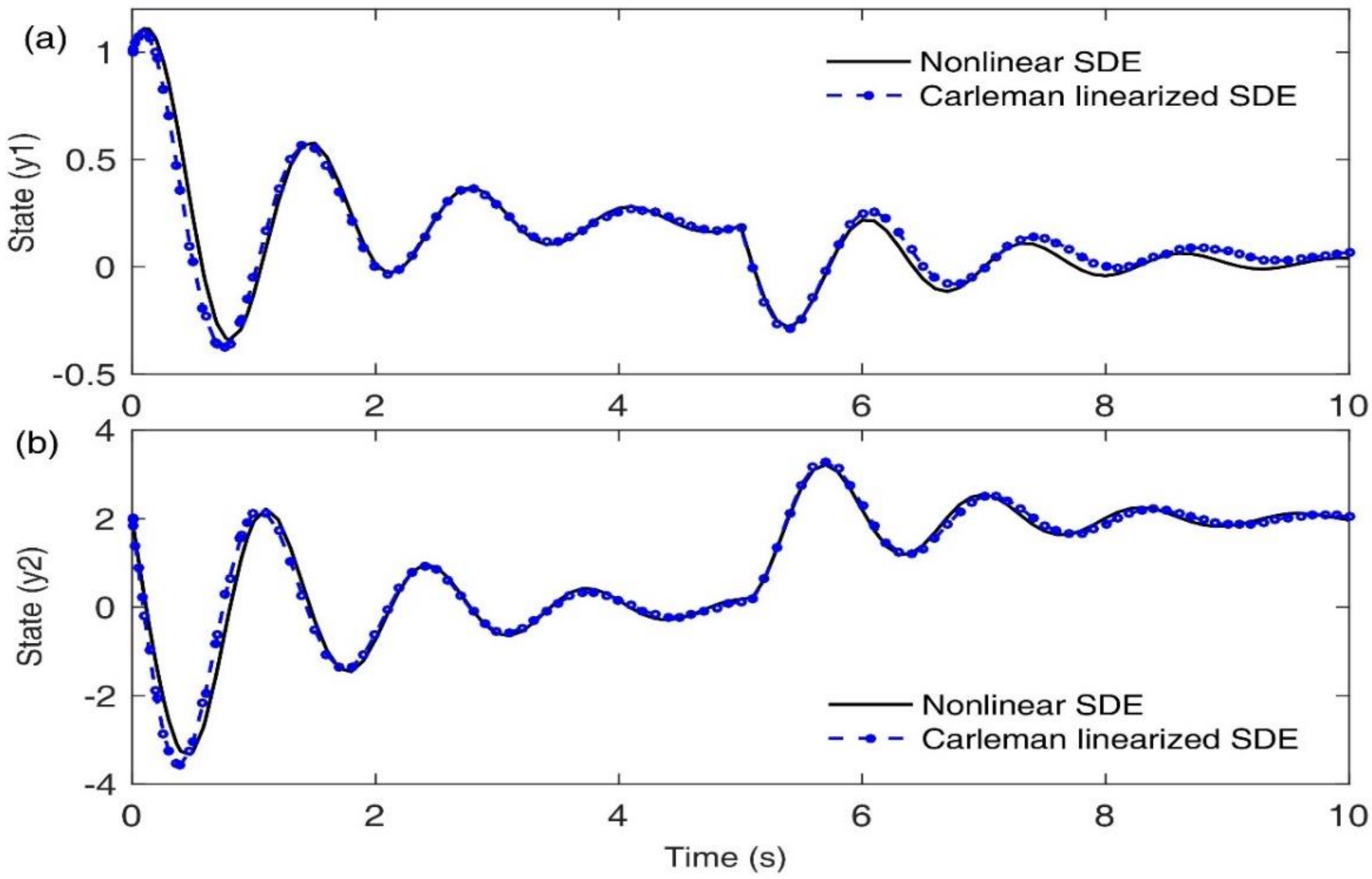

**Fig. 3:** A comparison between the machine swing SDE and its Carleman linearized SDE.

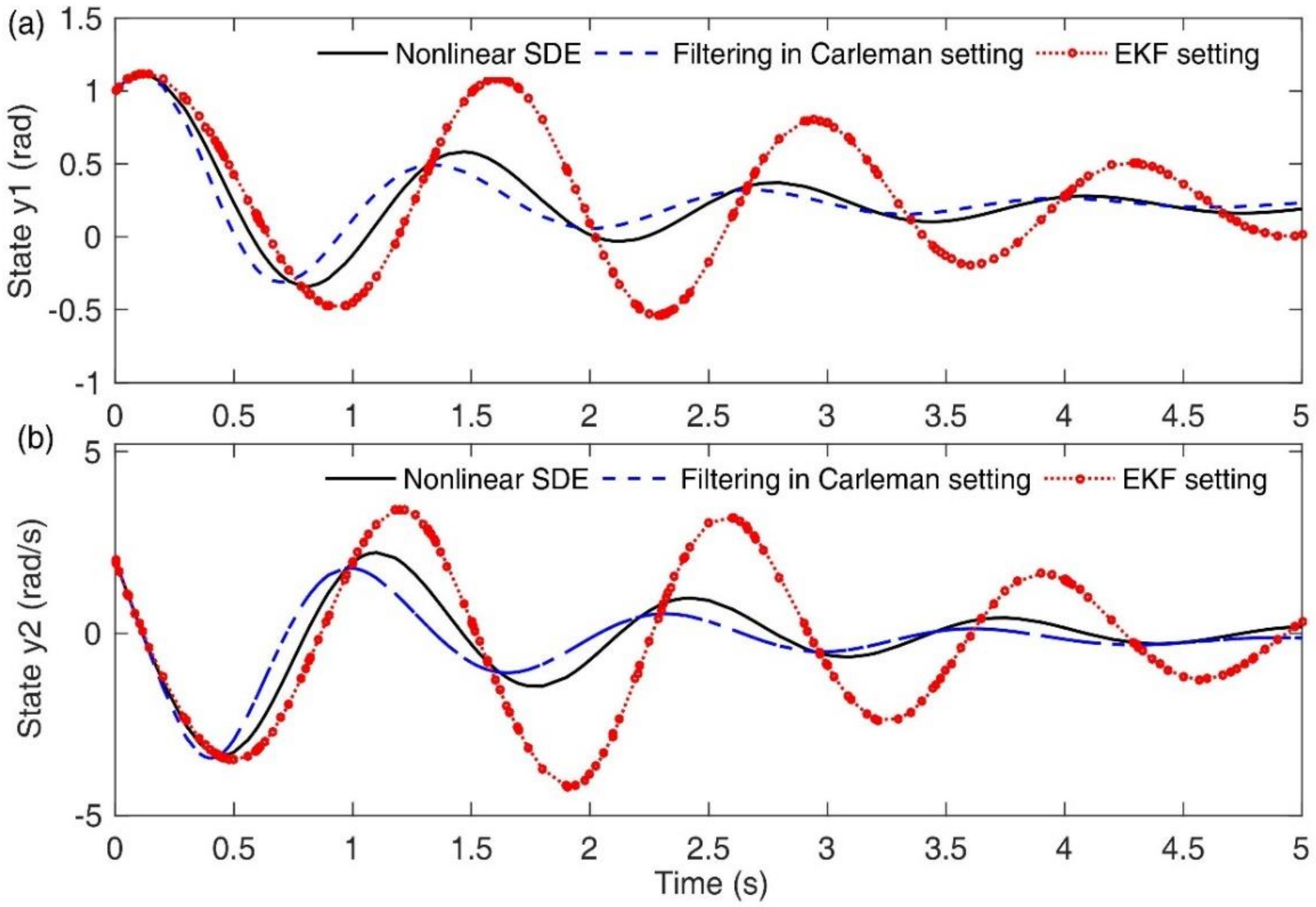

**Fig. 4:** (a) conditional mean trajectories for the state $y_1$, (b) conditional mean trajectories for the state $y_2$.

Fig. 3 reveals that the Carleman linearized trajectories are associated with the states $y_1$ and $y_2$ tracks the nonlinear state trajectories. The closeness of the linearized trajectories to that of the concerning noisy exact nonlinear trajectories reveals that the Carleman linearization of the order three for the machine swing dynamics captures the sinusoidal nonlinearity quite well.

Now, we demonstrate the results by utilizing a filtering set of equations (16)-(17) for the Carleman linearized machine swing SDE. Fig. 4 contrasts three trajectories. The first trajectory is concerning the true state trajectory (solid lines), an evolution of the nonlinear exact stochastic differential equation. The second trajectory is the Carleman linearized filtered state trajectory (dashed lines) resulting from the conditional mean equations (16a)-(16i). The third trajectory is the EKF filtered estimate state trajectory (marked dot-dot lines) resulting from the conditional mean equations (19a)-(19b).

Note that the Carleman linearized mean trajectories are generated by the 'algorithmic filtering procedure'. Fig. 4 reveals that the filtered state trajectories in the Carleman setting show 'the greater closeness with the actual associated state trajectories as well as the boundedness property'. The Carleman linearized filtered state trajectories are closer to the nonlinear SDE state trajectories in contrast to the EKF trajectories.

Filtering in the Carleman setting can be adjudged using the notion of the absolute filtering error via simulations, see Figs. (5a)-(5b). The absolute filtering error evolves for the Carleman setting as well as the EKF setting that shows the well-behaved increase and decrease. Importantly, the error is less for the Carleman linearized case and becomes quite larger for the EKF case. Table 2 shows a contrast for both the states. That is indicative of the Carleman framework for the filtering is superior to the EKF setting.

**Table 2:** A maximum absolute error comparison for the first set of parameters

| Filtering method | state $y_1$ max. absolute error | state $y_2$ max. absolute error |
|---|---|---|
| EKF setting | 0.65 | 0.34 |
| Proposed Carleman linearized setting | 0.22 | 0.10 |

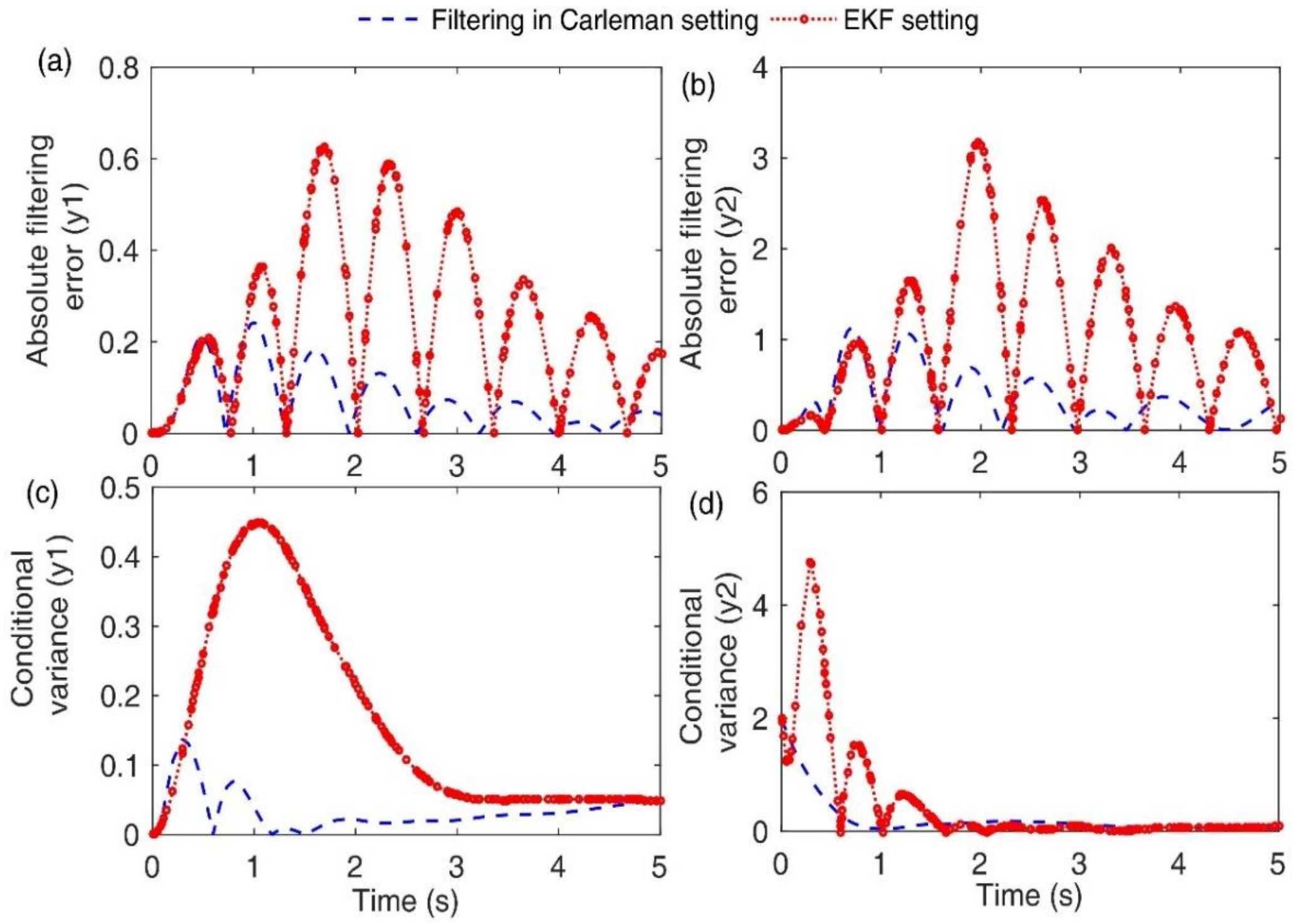

**Fig. 5:** Comparison of absolute filtering errors and conditional variances trajectories for the first set of parameters

Furthermore, Figs. 5(c) and 5(d) show an illustration of the conditional variance evolutions for both states. The illustration is made for two frameworks, the Carleman as well as the EKF. Since the expectation of the square of the filtering error is the conditional variance. The conditional variance is the positive semi-definite as well as the expectation of the non-negative stochastic processes is always non-negative. This implies the conditional variance associated with the Carleman setting is considerably less in contrast to the EKF. That is illustrated graphically in Figs. (5c)-(5d) as well.

### *4.2 Second Set*

For gaining insights into the usefulness of filtering in two settings, we exploit the second set of parameters. The second set is different from the first set in the sense of the mechanical design of the machine. The damping of the machine in the second set is larger than the first set. The associated initial conditions are $y_1(0)=1.5,\ y_2(0)=1,\ P_{y_1}(0)=0.1,\ P_{y_2}(0)=2$ and $P_{y_1 y_2}(0)=0.$ The second set of parameters is intended to achieve the numerical simulation of a relatively larger damping machine.

**Table 3:** The second set of machine swing dynamics parameters

| Parameters | $D$ | $E_t$ | $V$ | $X$ | $\sigma_1$ | $\sigma_2$ | $P_m$ | $M$ | $\varphi_n$ |
|---|---|---|---|---|---|---|---|---|---|
| Values | 0.7 | 1.2 | 1 | 0.25 | 0.08 | 0.06 | 1 | 0.2 | 100 |
| Units | $pu/rad/s$ | $pu$ | $pu$ | $pu$ | - | - | $pu$ | - | - |

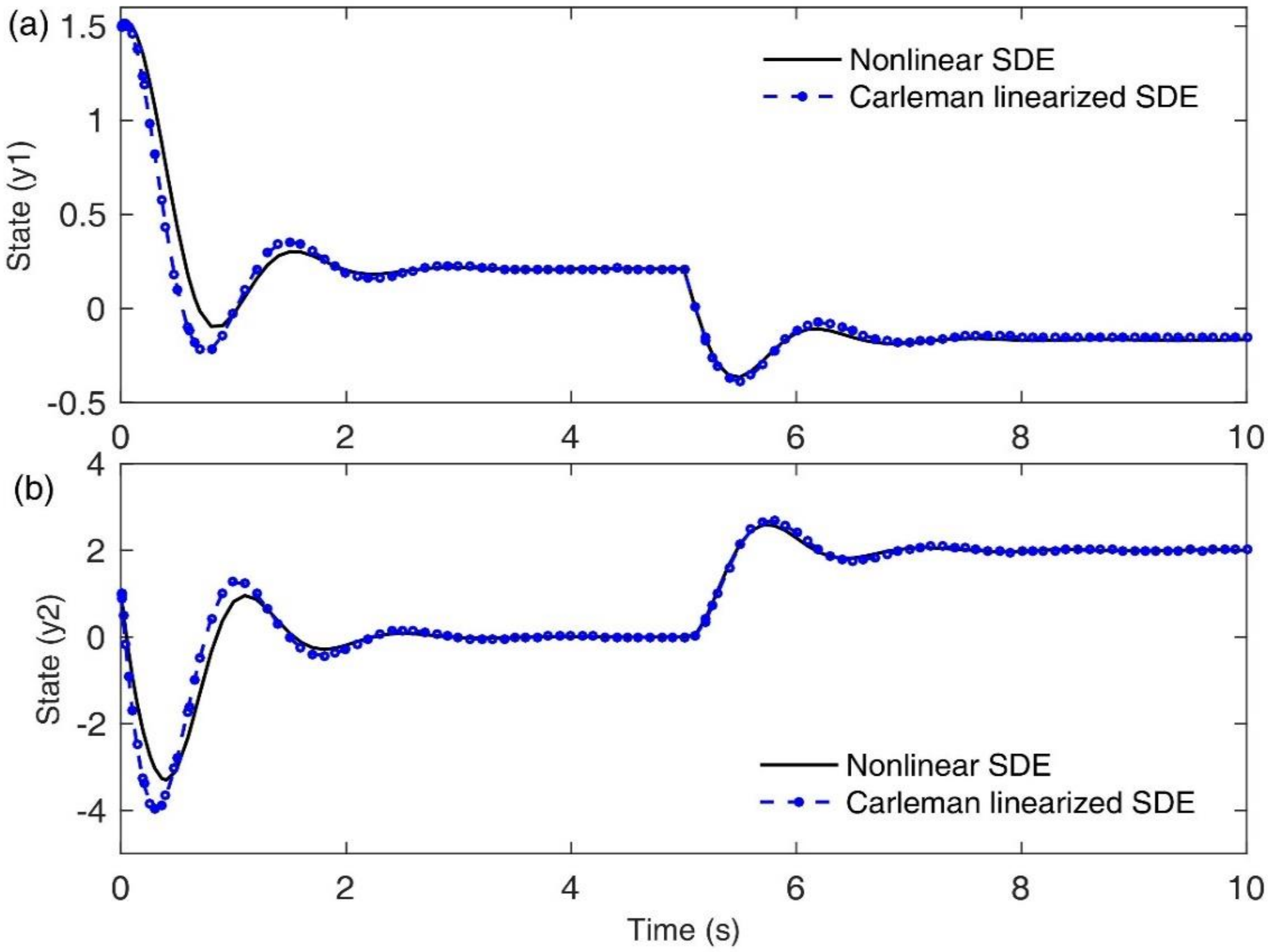

**Fig. 6:** A comparison between a machine swing SDE and its Carleman linearized SDE

The procedure to achieve filtering via numerical simulations for the 'second set' (Table 3) is the same as that of the first set. A step-change in the rotor angle $\delta_t$ is introduced at the time interval of 5 sec, via injecting a disturbance of the absolute magnitude 2 $rad.$ Simulation results for a comparison between the Carleman linearized state stochastic evolutions and the concerning exact stochastic evolutions are illustrated in Fig. 6.

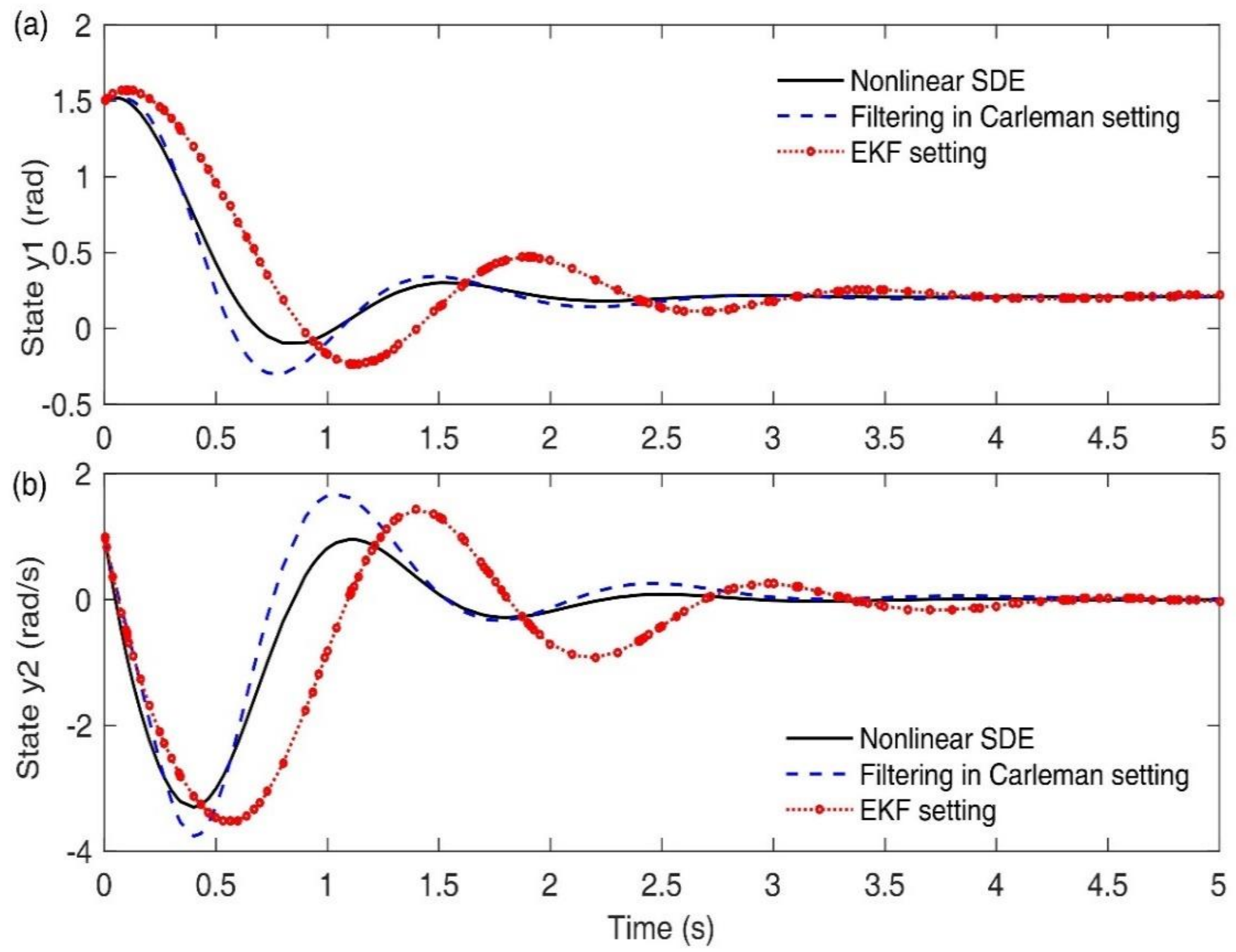

**Fig. 7:** (a) conditional mean trajectory for the state $y_1$, (b) conditional mean trajectory for the state $y_2$.

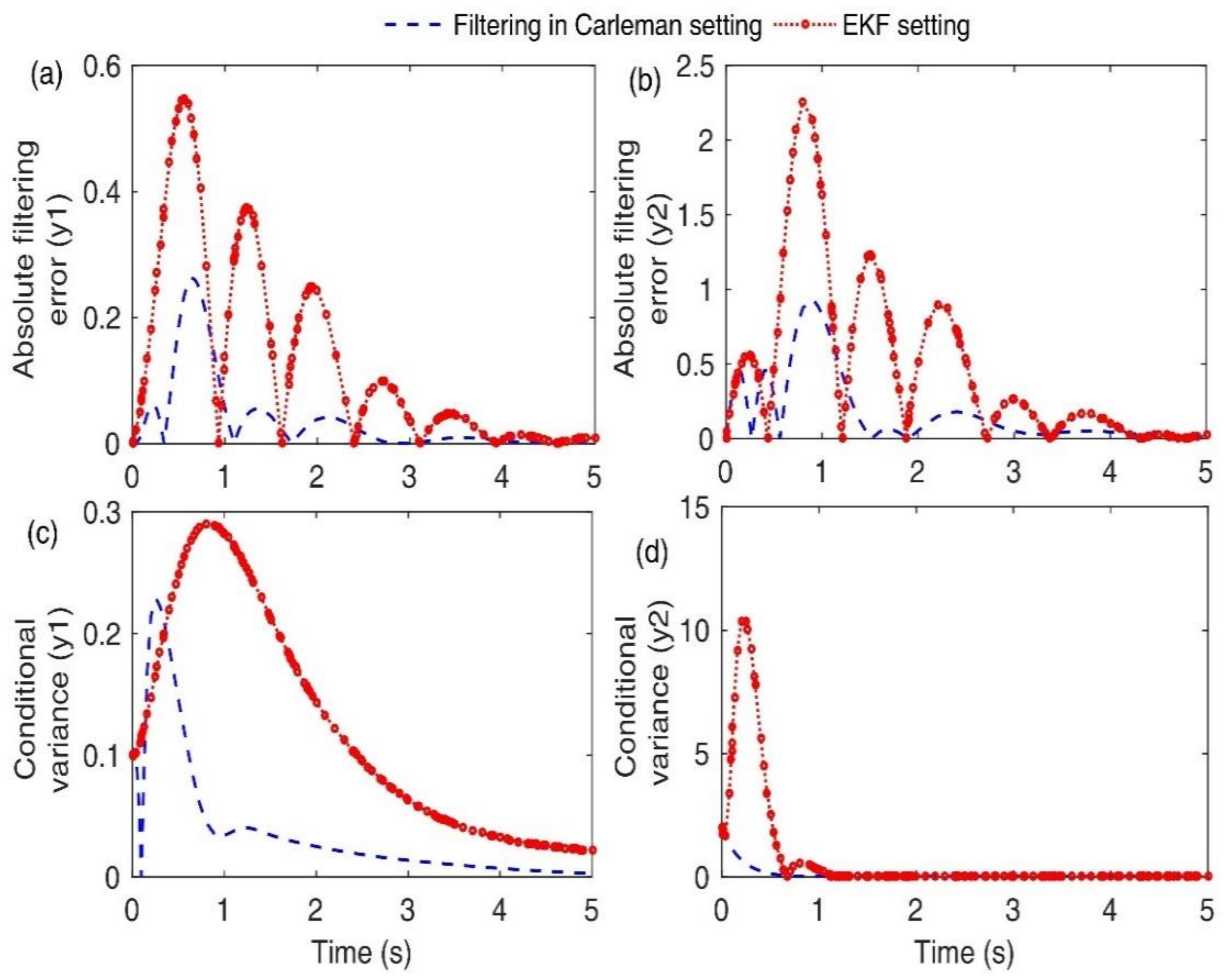

**Fig. 8:** Comparison of absolute filtering errors and conditional variances trajectories for the second set of parameters

Figure (7) unfolds the machine swing stochastically perturbed dynamics as well as the most probable dynamics stemming from nonlinear filtering techniques of stochastic processes. The contribution to the evolution of the state vector comes from initial conditions as well as random inputs. The most probable dynamics is a consequence of the action of the conditional expectation of the stochastic state. After a long-time of evolution, figure (7) shows the states do not evolve leading to the equilibrium point. Furthermore, the Carleman linearized filtered estimates track more closely the exact state in contrast to the EKF state.

Since this paper is chiefly intended to suggest the usefulness of the Carleman framework for the filtering problems. Similar to the first set of parameters, we exploit the second set of parameters and adopt the notion of the absolute filtering error from the perspectives of the Carleman framework as well as EKF. The second set of parameters in combination with the notion of the absolute filtering error infers the refinement in the estimates attributed to the Carleman framework. Fig. 8 reveals greater usefulness of the Carleman setting in contrast to the EKF.

# 5. Conclusion

The novelty of the paper is to achieve the Carleman linearization of a Markovian nonlinear stochastic differential system. To illustrate the idea of the Carleman linearization of Markov processes and achieving their filtering, this paper exploits the nonlinear machine swing dynamics as well.

Interestingly, the nonlinear machine swing dynamics is a Markovian dynamic. The usefulness of the swing equation filtering can be ascertained from the fact that extracting machine swing states from the on-line machine swing 'noisy' data. It is worth to mention that the swing equation filtering in the Carleman framework is refined and sharper in contrast to the benchmark EKF estimation framework. Extracting inaccurate states may lead to poor decisions about power system operations and leading to severe consequences. Thus, this paper will be of interest to nonlinear dynamists aspiring for refinements in algorithms for extracting 'nonlinear system states from the noisy data', including on-line machine swing noisy data.

The method of this paper is quite general, this paper recommends the usefulness of the Carleman framework for nonlinear filtering problems that are relatively very less researched.

# Appendix:

Consider the vector Gaussian stochastic process $Y=\{y_t,\Im_t,0\le t<\infty\}$, where $y_t\in R^n$. Furthermore $y_t$ is normal with mean $\mu_{y_t}$ and the variance $P_{y_t}$. Here we wish to construct the characteristic function and then weave the higher-order moment equations. The higher-order moment equation arises from practical problems involving higher-order statistics. That introduces mathematical subtleties. Here, we introduce the idea of replacing the higher-order moment with lower-order moments. The appendix does that. Now construct the compensated vector Gaussian process $\tilde{y}_t=y_t-\hat{y}_t$ where the compensated process is zero mean and the variance $P_{y_t}$. Thus, the concerning Gaussian characteristic function is

$$E\exp(s^T\tilde{y}_t)=\exp\left(\frac{1}{2}s^T P_{y_t}s\right)=\exp\left(\frac{1}{2}\sum_{i,j}s_i s_j P_{ij}\right)=\sum_m\frac{\left(\frac{1}{2}\sum_{i,j}s_i s_j P_{ij}\right)^m}{m!}$$

$$=\sum_m\frac{\left(\frac{1}{2}\sum_i s_i^2 P_{ii}+\sum_{i_1<i_2}s_i s_j P_{ij}\right)^m}{m!}. \tag{A.1}$$

From the generating functions,

$$\sum_m\sum_{r_1+r_2+\cdots+r_N=m}E\tilde{y}_1^{r_1}\tilde{y}_2^{r_2}\cdots\tilde{y}_N^{r_N}\frac{s_1^{r_1}}{r_1!}\frac{s_2^{r_2}}{r_2!}\cdots\frac{s_N^{r_N}}{r_N!}=E\exp(s^T\tilde{y}_t). \tag{A.2}$$

From equations (A.1) and (A.2), we have

$$\sum_{m}\sum_{r_1+r_2+\cdots+r_N=m} E\tilde{y}_1^{r_1}\tilde{y}_2^{r_2}\cdots\tilde{y}_N^{r_N}\frac{s_1^{r_1}}{r_1!}\frac{s_2^{r_2}}{r_2!}\cdots\frac{s_N^{r_N}}{r_N!}=\sum_{p}\frac{\left(\frac{1}{2}\sum_{i}s_i^2P_{ii}+\sum_{i_1<i_2}s_is_jP_{ij}\right)^p}{p!}.$$

Rearranging the above terms, we have

$$\sum_{m\text{ is even}}\sum_{r_1+r_2+\cdots+r_N=m} E\tilde{y}_1^{r_1}\tilde{y}_2^{r_2}\cdots\tilde{y}_N^{r_N}\frac{s_1^{r_1}}{r_1!}\frac{s_2^{r_2}}{r_2!}\cdots\frac{s_N^{r_N}}{r_N!}+\sum_{m\text{ is odd}}\sum_{r_1+r_2+\cdots+r_N=m} E\tilde{y}_1^{r_1}\tilde{y}_2^{r_2}\cdots\tilde{y}_N^{r_N}\frac{s_1^{r_1}}{r_1!}\frac{s_2^{r_2}}{r_2!}\cdots\frac{s_N^{r_N}}{r_N!}$$

$$=\sum_{p}\frac{\left(\frac{1}{2}\sum_{i}s_i^2P_{ii}+\sum_{i_1<i_2}s_is_jP_{ij}\right)^p}{p!}.$$

The above suggests that the second term of the left-hand side vanishes.

$$\sum_{m\text{ is even}}\sum_{r_1+r_2+\cdots+r_N=m} E\tilde{y}_1^{r_1}\tilde{y}_2^{r_2}\cdots\tilde{y}_N^{r_N}\frac{s_1^{r_1}}{r_1!}\frac{s_2^{r_2}}{r_2!}\cdots\frac{s_N^{r_N}}{r_N!}=\sum_{p}\frac{\left(\frac{1}{2}\sum_{i}s_i^2P_{ii}+\sum_{i_1<i_2}s_is_jP_{ij}\right)^p}{p!}.$$

Suppose $m=q$, where $m$ varies and $q$ is fixed and even, then $p=\frac{q}{2}$ and

$$\sum_{r_1+r_2+\cdots+r_N=q} E\tilde{y}_1^{r_1}\tilde{y}_2^{r_2}\cdots\tilde{y}_N^{r_N}\frac{s_1^{r_1}}{r_1!}\frac{s_2^{r_2}}{r_2!}\cdots\frac{s_N^{r_N}}{r_N!}=\frac{\left(\frac{1}{2}\sum_{i}s_i^2P_{ii}+\sum_{i_1<i_2}s_is_jP_{ij}\right)^{\frac{q}{2}}}{\frac{q}{2}!}. \quad \text{(A.3)}$$

Here, three cases arise, $q>N,\ q=N,\ q<N.$ Here, we restrict our discussions to $q=N$ and we pose a question to compute the sixth-order moment explicitly in terms of the second-order under the Gaussian statistics, i.e.

$$E(y_1(t)-\hat{y}_1(t))(y_2(t)-\hat{y}_2(t))(y_3(t)-\hat{y}_3(t))(y_4(t)-\hat{y}_4(t))(y_5(t)-\hat{y}_5(t))(y_6(t)-\hat{y}_6(t))=?$$

To answer the above, set

$$q=6,\ N=6,\ r_1+r_2+\ldots+r_6=6,\ r_1=1,\ r_2=1,\ldots,r_6=1,$$

and invoke the condition

$$\bigcap_k\left\{(i_k,i_{k+1})\middle|1\le i_k\le 6, 1\le i_{k+1}\le 6,\ i_k<i_{k+1}\right\}=\phi.$$

That can be achieved using the following: first, construct the set $\left\{(i_1,i_2)\middle|1\le i_1\le 6,1\le i_2\le 6,1\le i_1<i_2\le 6\right\}$. The set is a product space and each element of the product space is a set of two tuples in increasing order. Now, we construct a product space with the property $\bigcap_k\left\{(i_k,i_{k+1})\middle|1\le i_k\le 6, 1\le i_{k+1}\le 6,\ i_k<i_{k+1}\right\}=\phi,$ where $k$ is odd valued with $1\le k\le 6.$ After a result of these, we equate the terms of the left-hand and right-hand sides of equation (A.3), associated with the term $s_1s_2s_3s_4s_5s_6$. Thus, we have

$$E\prod_{1\le i\le 6}(y_i(t)-\hat{y}_i(t))=P_{y_1y_2}P_{y_3y_4}P_{y_5y_6}+P_{y_1y_2}P_{y_3y_5}P_{y_4y_6}+P_{y_1y_2}P_{y_3y_6}P_{y_4y_5}+P_{y_1y_3}P_{y_2y_4}P_{y_5y_6}$$
$$+P_{y_1y_3}P_{y_2y_5}P_{y_4y_6}+P_{y_1y_3}P_{y_2y_6}P_{y_4y_5}+P_{y_1y_4}P_{y_2y_3}P_{y_5y_6}+P_{y_1y_4}P_{y_2y_5}P_{y_3y_6}$$
$$+P_{y_1y_4}P_{y_2y_6}P_{y_3y_5}+P_{y_1y_5}P_{y_2y_3}P_{y_4y_6}+P_{y_1y_5}P_{y_2y_4}P_{y_3y_6}+P_{y_1y_5}P_{y_2y_6}P_{y_3y_4}$$

$$+P_{y_1y_6}P_{y_2y_3}P_{y_4y_5}+P_{y_1y_6}P_{y_2y_4}P_{y_3y_5}+P_{y_1y_6}P_{y_2y_5}P_{y_3y_4}. \tag{A.4}$$

Consider the scalar Gaussian stochastic process, the above equation (A.4) is simplified to

$$E(y_t-\hat{y}_t)^6=15P_{y_t}^3. \tag{A.5}$$

Furthermore, we wish to calculate $Ey_t{}^6$. That can be calculated using the binomial coefficient as well as lower-order Gaussian statistics. The relation holds:

$$E(y_t-\hat{y}_t)^6=E\sum_{r=0}^{6}(-1)^{6-r}\binom{6}{r}y_t^r\hat{y}_t{}^{6-r}=\sum_{r=0}^{6}(-1)^{6-r}\binom{6}{r}Ey_t^r\,\hat{y}_t{}^{6-r}$$

$$={}^6C_0\hat{y}_t{}^6-{}^6C_1Ey_t\hat{y}_t{}^5+{}^6C_2Ey_t^2\hat{y}_t{}^4-{}^6C_3Ey_t^3\hat{y}_t{}^3+{}^6C_4Ey_t^4\hat{y}_t{}^2-{}^6C_5Ey_t^5\hat{y}_t+{}^6C_6Ey_t^6.$$

After combining equation (A.5) with the above, we have

$${}^6C_0\hat{y}_t{}^6-{}^6C_1Ey_t\hat{y}_t{}^5+{}^6C_2Ey_t^2\hat{y}_t{}^4-{}^6C_3Ey_t^3\hat{y}_t{}^3+{}^6C_4Ey_t^4\hat{y}_t{}^2-{}^6C_5Ey_t^5\hat{y}_t+{}^6C_6Ey_t^6=15P_{y_t}^3.$$

Thus,

$$Ey_t^6=15P_{y_t}^3-{}^6C_0\hat{y}_t{}^6+{}^6C_1Ey_t\hat{y}_t{}^5-{}^6C_2Ey_t^2\hat{y}_t{}^4+{}^6C_3Ey_t^3\hat{y}_t{}^3-{}^6C_4Ey_t^4\hat{y}_t{}^2+{}^6C_5Ey_t^5\hat{y}_t. \tag{A.6}$$

Using the binomial coefficient identities and expectation operator, i.e.

$$Ey_t^3=\left(3\hat{y}_tP_{y_t}+\hat{y}_t^3\right),\quad Ey_t^4=\left(3P_{y_t}^2+6\hat{y}_t^2P_{y_t}+\hat{y}_t^4\right),\quad Ey_t^5=\left(15\hat{y}_tP_{y_t}^2+10\hat{y}_t^3P_{y_t}+\hat{y}_t^5\right),$$

as well as embedding them into equation (A.6), we have

$$Ey_t^6=\left(15P_{y_t}^3+45\hat{y}_t^2P_{y_t}^2+15\hat{y}_t^4P_{y_t}+\hat{y}_t{}^6\right).$$